\documentclass{IEEEtran}

\IEEEoverridecommandlockouts

\usepackage{amsmath,amssymb,amsfonts}
\usepackage[capitalise]{cleveref}
\usepackage[]{footmisc}
\usepackage{xcolor}
\usepackage{algorithm}
\usepackage{algorithmic}
\usepackage[algo2e]{algorithm2e} 
\usepackage{graphicx}
\usepackage{textcomp}
\usepackage{psfrag}
\usepackage{epsfig}
\usepackage{subcaption}
\usepackage{float}
\usepackage{xcolor}
\usepackage{nicematrix}
\usepackage[numbers,sort&compress]{natbib}
\usepackage{tabularx}
\usepackage{booktabs}
\usepackage{diagbox}
\usepackage{multicol}
\newtheorem{mydef}{Definition}
\newtheorem{mypro}{Proposition}

\newtheorem{myrem}{Remark}
\newtheorem{mylem}{Lemma}
\newtheorem{mythr}{Theorem}

\begin{document}

\title{\LARGE Interior-Point-based $H_2$ Controller Synthesis for Compartmental Systems}

\author{Zhaohua~Yang,~Nachuan~Yang,~Pengyu~Wang,~Haishan~Zhang,~Xiayan~Xu,~and~Ling~Shi,~\IEEEmembership{Fellow,~IEEE}% \prec-this % stops a space
	\thanks{\rm Z. Yang, N. Yang, P. Wang, H. Zhang, X. Xu, and L. Shi are with the Department of Electronic and Computer Engineering, Hong Kong University of Science and
		Technology, Clear Water Bay, Hong Kong SAR (email: zyangcr@connect.ust.hk; nc.yang@connect.ust.hk; pwangat@connect.ust.hk; hzhangdo@connect.ust.hk; xxucj@connect.ust.hk; eesling@ust.hk).}
	\thanks{}
	 %\thanks{\rm Karl H. Johansson is with the Division of Decision and Control Systems, School of Electrical Engineering and Computer Science, KTH Royal Institute of Technology, 10044, Stockholm, Sweden (email: kallej@kth.se).}
  }
\maketitle

\begin{abstract}
This paper focuses on the optimal $H_2$ controller design for compartmental systems, with the aim of enhancing system robustness while maintaining the law of mass conservation. Through a novel problem transformation, we establish that the original problem is equivalent to a new non-convex optimization problem with a closed polyhedron constraint. Existing works have developed various first-order methods to tackle inequality constraints. 
% However, the performance of the first-order method is limited in terms of convergence speed and precision, restricting its potential in practical applications. Therefore, developing a novel algorithm with fast speed and high precision is critical. 
However, they often lack convergence guarantees in non-convex scenarios, thereby reducing their reliability in practical applications. Consequently, there is a critical need to develop new and efficient algorithms with convergence guarantees.
In this paper, we reformulate the problem using log-barrier functions and introduce two separate approaches with convergence guarantee to address the problem: first-order 
 interior point method (FIPM) and second-order interior point method (SIPM). Additionally, we propose an initialization method to guarantee the interior property of initial values. Finally, we compare FIPM and SIPM through a room temperature control example and show their pros and cons. They are also compared to the existing alternating direction method of multipliers (ADMM) across different system scales.
\vspace{3pt}
\end{abstract}


\IEEEpeerreviewmaketitle
\section{Introduction}
Compartmental systems refer to systems in which units, called compartments, interact with each other, and meanwhile obey the law of mass conservation. Such systems were first introduced in \cite{jacquez1972compartmental}. A more general introduction can be found in \cite{haddad2010nonnegative}. In many industrial systems with no interaction with the environment, the descriptor variables are not only non-negative but also their sum is non-increasing. For example, without the external supply, the total volume of material satisfies the law of mass conservation in a discrete-time compartmental flow model \cite{chang1977discrete}. Besides, the air traffic network can be divided into many air centers with conservative inflows and outflows\cite{deng2023event}. Similarly, the network that models the vehicles in the highway can also be modeled to be compartmental\cite{coogan2015compartmental}. Compartmental systems are a special kind of positive systems that can easily be implemented in reality. A classical tutorial on positive systems can be found in \cite{rantzer2018tutorial}.

$H_2$ and linear quadratic
regulator (LQR) performance are two important performance criteria in the control field. LQR performance characterizes the system cost by a quadratic function, while $H_2$ performance reflects the energy (variance) of the system output under unit-variance white-noise input. It can be shown that these two performance metrics are mathematically equivalent. Recently, the $H_2$ and LQR control with the positive constraint have attracted many researchers in the control community. Deaecto and Geromel \cite{Deaecto} proposed an iterative linear matrix inequality (ILMI) approach that can generate a solution sequence with non-increasing $H_2$ performance, thus guaranteeing convergence to a sub-optimal point. Furthermore, Ebihara et al. \cite{ebihara2019h_} derived the upper bound and the lower bound of the optimal $H_2$ performance using semidefinite programming (SDP). Recently, Yang et al.\cite{yang2022linear} proposed a projection-based LQR control of positive systems and showed it outperforms the methods in \cite{Deaecto}. As a special kind of positive systems, some results have been obtained on the stabilization of compartmental systems recently. Valcher and Zorzan \cite{valcher2018state} and Leenheer and Aeyels \cite{de2001stabilization} thoroughly studied the stabilization of compartmental systems and provided feasible approaches to find the stabilized controller under different circumstances. Different from previous works, we are interested in reducing the effects of disturbances on system outputs for a compartmental system that preserves physical property. In our recent work \cite{yang2024h2}, the optimal $H_2$ control of continuous-time compartmental systems was discussed, but the algorithm alternating direction method of multipliers (ADMM) is first-order without convergence guarantee. Meanwhile, the discrete-time case is still missing, and more efficient computation methods, including higher-order algorithms, remain to be developed. This motivates our research in this paper.

Extensive research has been conducted on the LQR and $H_2$ controller synthesis subject to structural constraints. 
Lin et al. \cite{lin2011augmented} proposed an augmented Lagrangian method to design an optimal $H_2$ controller with a given pattern. 
 In addition. Chanekar et al. \cite{chanekar2017optimal} considered adding a lower bound and an upper bound on the controller. Furthermore, Wu \cite{wu2021design} considered general structural constraints, i.e., the constraints are mixed with equality, inequality, and scaling.
% Furthermore, Wu \cite{wu2022structured} considered the mixed $H_2$/$H_\infty$ control where the $H_\infty$ is kept within a threshold. 
However, these works only relied on the first-order method and did not consider the Hessian information. Fatkhullin and Polyak \cite{polyak} proposed a Hessian operator acting in a specific direction. However, this method can only be applied to accelerate the line search. The main algorithm is still a first-order method. Recently, Cheng et al. \cite{cheng2022second} proposed a second-order method in the continuous case, but only equality constraints are considered. To the best of our knowledge, no existing work considered the second-order method on the design of controllers subject to compartmental constraints, a kind of linear inequality constraint that cannot be tackled by the aforementioned methods.

In this paper, our objective is to efficiently and reliably design the optimal $H_2$ controller for compartmental systems. This contribution is partitioned into multiple steps to achieve our goal. First, we propose a novel problem transformation technique and establish an equivalent optimization problem with a closed polyhedron feasible region. Second, we characterize
the optimality conditions of the reformulated problem via the Karush-Kuhn-Tucker (KKT) conditions and define stationary points. Third, we introduce, for the first time, the interior point method with the log-barrier function and propose a first-order 
 interior point method (FIPM) and a novel second-order 
 interior point method (SIPM). We establish their convergence to a stationary point of the original problem. Fourth, We develop a procedure to determine a strictly feasible controller for initialization. Finally, we conduct comprehensive numerical simulations to compare the performance of FIPM, SIPM, and ADMM \cite{yang2024h2} across various system dimensions.

This paper is organized as follows. In \cref{Preliminaries}, some preliminaries on compartmental systems and $H_2$ performance are provided and the problem
is formulated. In \cref{Basic Results}, we conduct the problem transformation and present the first-order optimality conditions. In \cref{First-Order Method}, we introduce the log-barrier term and propose the FIPM with guaranteed convergence. In \cref{Second-Order Method}, we provide detailed derivations for the Hessian matrix of the $H_2$ objective function and the log-barrier term, propose a Hessian modification method to guarantee descending and present the SIPM with guaranteed convergence. In \cref{Discussions on the initialization of K0}, we discuss the initialization of a strictly feasible controller. In \cref{Simulations}, several numerical simulations on room temperature control are provided to compare the performance of FIPM, SIPM and ADMM \cite{yang2024h2} across different system scales. The paper is concluded in \cref{Conclusion}.

\textit{Notations}:
{The notation} $R^{m\times n}$ denotes the set of all real matrices of size $m\times n$, and $R^{m\times n}_+$ denotes the set of real matrices with non-negative entries of size $m\times n$. For a scalar $t$, $\{t^{(k)}\}$ denotes a scalar sequence. For a matrix $X$, $X^T$ denotes its transpose, $X^H$ denotes its conjugate transpose, $X^k$ denotes its $k^{th}$ power, $\{X^{(k)}\}$ denotes a matrix sequence, $Tr(X)$ denotes its trace, $X_{ij}$ denotes its element at $i$-th row and $j$-th column. The notation $\partial$ denotes the boundary of a set. The identity matrix and the zero matrix are denoted by $I$ and 0, with dimensions labeled using subscripts if necessary. 
% When each entry of a matrix $X$ is positive (resp. non-negative), we write $X>0$ (resp. $X\ge 0$). When a matrix $X$ is positive semidefinite(resp. definite), we write $X\succeq 0$ (resp. $X\succ 0$). 
The notation $X>0$ (or $X\ge 0$) means each entry of matrix $X$ is positive (or non-negative). The notation $X\succ 0$ (or $X\succeq 0$) means matrix $X$ is positive definite (or semidefinite). We use $\mathbf{0}$ and $\mathbf{1}$ to denote the all-zeros column vectors and all-ones column vectors with compatible dimensions. $\mathrm{diag}(a_1,a_2,\ldots,a_n)$ represents a diagonal matrix with $a_i, \forall{i}=1,2,\ldots,n$ on its diagonal. Similarly, $\mathrm{blkdiag}(A_1,A_2,\ldots,A_n)$ represents putting $A_i, \forall{i}=1,2,\ldots,n$ in a block diagonal way. $\left[A_1|A_2|\ldots|A_n\right]$ represents organizing $A_i, \forall{i}=1,2,\ldots,n$ in one row. The operator $\otimes$ represents the Kronecker product. $(\cdot)^*$ represents the adjoint operator. $\circ$ represents the Hadamard power (elementwise power) \cite{reams1999hadamard}. 
The operator $\mathrm{vec}(\cdot)$ represents the vectorization operation that expands a matrix by column into a column vector. $\mathrm{mat}(\cdot)$ represents the matrixization operation that reorganizes a column vector into a matrix with appropriate dimensions by column.
% The operator $\mathrm{vec}(\cdot)$ represents  expanding a matrix by column into a column vector. $\mathrm{mat}(\cdot)$ represents reorganizing a column vector into a matrix with appropriate dimensions by column.

\section{Preliminaries}\label{Preliminaries}
\subsection{Compartmental Systems}
We can use a non-negative matrix $A\in R^{n\times n}$ to describe a digraph $\mathcal D(A) = \{\mathcal V, \mathcal E\}$, where $\mathcal V = \{1,...,n\}$ is the set of vertices and $\mathcal E$ is the set of edges. There is an arc $(i,j)\in \mathcal E$ from $i$ to $j$ if and only if $A_{ji}>0$, where $A_{ji}$ is called the weight of the arc. A sequence of nodes $i_1 \rightarrow i_2 \rightarrow \cdots \rightarrow i_n$ is called a path if $(i_1,i_2)\in \mathcal E,\ldots,(i_{n-1},i_n) \in \mathcal E$. We say $j$ is accessible from $i$ if there is a path from $i$ to $j$, or equivalently, $A^l_{ji}>0$ for some $l$. Two different nodes $i$ and $j$ communicate if each of them is accessible from the other. Therefore, we can divide $\mathcal V$ into communicating classes, say $C_1, C_2 \ldots C_k$. 
 If there exists $g \in C_i$ and $h \in C_j$ such that $g$ accesses $h$, then we say $C_i$ accesses $C_j$. Each $C_i$ accesses itself trivially. The digraph $\mathcal D(A)$ is strongly connected if every two nodes communicate with the other. Equivalently, there is only one communicating class $\mathcal V$. The digraph $\mathcal D(A)$ is strongly connected if and only if $A$ is irreducible \citep[Definition 6.2.21, Definition 6.2.22, Theorem 6.2.24]{horn2012matrix}.

 \begin{mydef}\label{schur definition}
     A matrix $X$ is a Schur (stable) matrix if all its eigenvalues exist within an open unit circle. 
 \end{mydef}
A non-negative matrix, where each column sums to no more than 1 (i.e., $\mathbf{1}^TA\le \mathbf{1}^T$), is
called a compartmental matrix. This matrix is commonly utilized to describe flows between
compartments, illustrating an outflow from each node to the others. If the outflow is more than the inflow, it indicates a material loss to the environment. Mathematically, we call node $i$ an \textit{outflow node} if $\mathbf{1}^TAe_i<1$. The node $j$ is said to be \textit{outflow-connected} if there is a path from that node to an outflow node. Consequently, it is straightforward to deduce that a compartmental matrix is Schur if and only if every node is outflow-connected.

\begin{mydef}    
For a general linear time-invariant (LTI) system $x_{k+1} = Ax_k$, it is called a compartmental system if $A$ is a compartmental matrix.
% \begin{equation}
%     % % \begin{aligned}
%     x_{k+1} = Ax_k
%     % % % y_k &= Cx_k + Du_k\\
% % \end{aligned}
% \end{equation}
% where $A\in R^{n\times n}$ is a compartmental matrix.
\end{mydef}

% A compartmental system is illustrated in \cref{flow network}. The network is divided into two strongly connected communicating classes $C_1$ and $C_2$. Besides, $C_1$ accesses $C_2$ since there is a flow from node 4 to node 7. Further, node 5 is an outflow node. This compartmental system is Schur since every node is outflow-connected. 
%  \begin{figure}[t]
%     \centering
%     \includegraphics[scale=0.65]{flow network.png}
%     \caption{Illustration of the compartmental system}
%     \label{flow network}
% \end{figure}
\subsection{$H_2$ Performance}
Consider a discrete-time LTI system 
\begin{equation}
\label{basic system}
\begin{aligned}
    x_{k+1} &= Ax_k + Bu_k + Gd_k\\
    y_k &= Cx_k + Du_k\\
\end{aligned}
\end{equation}
where $x\in R^n$ is the state, $u\in R^m$ is the control, $d\in R^n$ is the exogenous disturbance, $y\in R^r$ is the sensor measurement. It is standard to assume that $(A, B)$ is stabilizable, $D^TC=0$ and $D^TD\succ0$ \cite{dullerud2013course}. The system adopts static state-feedback control law
\begin{equation}
\label{static state-feedback controller}
    u_k = -Kx_k
\end{equation}
where $K$ is the gain to be determined. Therefore, we can rewrite system (\ref{basic system}) into following form
\begin{equation}
    \label{rewritten system}
    \begin{aligned}
        x_{k+1} &= A_K x_k + Gd_k\\
        y_k &= C_K x_k
    \end{aligned}
\end{equation}
where $A_K=A-BK$ and $C_K=C-DK$. The closed-loop transfer function in $z$-domain from $d$ to $y$ is given by 
\begin{equation}
    T(z) = C_K(zI-A_K)^{-1}G
\end{equation}
It is known that the $H_2$ norm is bounded if the system is asymptotically stable. According to \cite{chen1995linear}, system \eqref{rewritten system} is asymptotically stable if and only if $A-BK$ is Schur. 
% According to \cite{chen1995linear}, system \eqref{rewritten system} is asymptotically stable if and only if all 
%  eigenvalues of $A-BK$ are within an open unit circle. 
 % To facilitate discussions, we introduce the definition of Schur matrix as follows.
 % \begin{mydef}\label{schur definition}
 %     A matrix $X$ is a Schur (stable) matrix if all its eigenvalues exist within an open unit circle. 
 % \end{mydef}
 
 With \cref{schur definition}, the implicit domain of $K$ is $K\in \mathcal{F}$, where $\mathcal{F}\triangleq\{K \mid A-BK~\text{is Schur}\}$. The $H_2$ norm for a stable transfer function $T(z)$ is defined by 
\begin{equation}
    \label{H2 norm}
    ||T(z)||_2 = \sqrt{\frac{1}{2\pi}\int^{\pi}_{-\pi}Tr(T(e^{-j\omega})^H T(e^{j\omega})) d\omega}
\end{equation}
By defining $J(K)\triangleq ||T(z)||^2_2$ as the objective function to be optimized, it is known \cite{peres1993h2} that we can extend the domain of $J(K)$ and redefine it as follows
\begin{equation}
        \label{J(K)}
        J(K) = \begin{cases}
            Tr(G^TXG), &K \in \mathcal{F}\\
            \infty, &\text{otherwise}
        \end{cases}
\end{equation}
where $X\succeq0$ is the unique solution to 
\begin{equation}
\label{X}
    A_K^TXA_K -X+C_K^TC_K=0
\end{equation}
\subsection{Problem Formulation}
In this paper, we are interested in the optimal static state-feedback controller synthesis under $H_2$ performance for compartmental systems. The problem is formulated below.

\textbf{Problem DHSCCS} (Design of $H_2$ Static State-Feedback Controller for Compartmental Systems): Consider the compartmental system in \eqref{basic system}, i.e., $A$ is compartmental, design a static state-feedback controller in \eqref{static state-feedback controller} such that the square of $H_2$ norm of the transfer function, as defined in \eqref{J(K)}, is minimized, and meanwhile the closed-loop system is both asymptotically stable and compartmental, i.e., $A_K$ is Schur and compartmental.
\begin{myrem}
    Intuitively, this problem is to find a $H_2$ static state-feedback controller that preserves the system's physical property, i.e., the flow conservation law of the compartmental systems, and meanwhile enhances the system's robustness against external disturbance. The main difficulty of this problem lies in three parts. First, the interior-point method requires a strictly feasible initial point. However, initializing a controller to ensure the system is both Schur and compartmental poses a significant challenge. Second, the $H_2$ performance is typically a non-convex function concerning $K$. Lastly, it is challenging for the system to simultaneously maintain Schur stability and compartmental properties. In this paper, we tackle these challenges from an optimization perspective using FIPM and SIPM. In numerical simulations, we will analyze their performance in terms of complexity, show their advantages and disadvantages, and compare them with ADMM \cite{yang2024h2}.
\end{myrem}

\section{Basic Results}\label{Basic Results}
This section presents fundamental results essential for subsequent advanced discussions. First, we transform the original problem into an optimization problem with polyhedron constraints. Second, we derive the gradient of $J(K)$, which is crucial for the first-order and the second-order methods. Finally, we depict the first-order optimal conditions via KKT conditions and define stationary points via these conditions.
\begin{mythr}
    Problem DHSCCS is equivalent to the following constrained optimization problem
    \begin{equation}
        \label{equivalent problem}
        \begin{aligned}
        &\mathop{\min}_{K} J(K)\\
        {s.t.}&\left\{\begin{array}{l} 
            A_K\geq 0,\\
            \mathbf{1}^TA_K \leq \mathbf{1}^T.
        \end{array}\right.
        \end{aligned}
        \end{equation}
\end{mythr}
\quad \textit{Proof}: Given $A_K$ must be a compartmental matrix, initially non-negative, the constraint $A_K\geq 0$ in \eqref{equivalent problem} is straightforward. When the compartmental matrix $A_K$ is Schur, ensuring every node is outflow connected, it follows that $\mathbf{1}^TA_K \leq \mathbf{1}^T$, with equality not holding. By defining $\mathcal{S} \triangleq \{K \mid A_K \text{ is compartmental and Schur}\}$ and the set of compartmental gains $\mathcal{C} \triangleq \{K \mid A_K\geq 0 \text{ and } \mathbf{1}^TA_K \leq \mathbf{1}^T\}$, it shows that $\mathcal{S} \subset \mathcal{C}$, $\mathcal{S} \subset \mathcal{F}$, and $\mathcal{C} - \mathcal{S} \subset \mathcal{\bar F}$, where $\mathcal{\bar F}$ denotes the complementary set of $\mathcal{F}$. Further note that as $K \rightarrow \partial\mathcal{F}$, $J(K) \rightarrow \infty$ due to the coercive property of $J(K)$ \cite{bu2019lqr}. Therefore, from an optimization perspective, we can relax $\mathcal{S}$ to $\mathcal{C}$ without loss of generality, and thus the \textbf{Problem DHSCCS} is equivalent to \eqref{equivalent problem}. The proof is complete. $\hfill \square$

The differentiability of $J(K)$ is a well-known fact, shown by Levine and Athans \cite{levine}. In what follows, we characterize the gradient of $J(K)$ in the discrete-time $H_2$ control case.
\begin{mylem}
    The gradient of $J(K)$ is given by 
    \begin{equation}
        \label{nabla J(K)}
        \nabla J(K) = 2(B^TX(A-BK)-D^TDK)Y
    \end{equation}
    where $Y\succeq0$ is the unique solution to 
    \begin{equation}
    \label{Y}
        A_KYA_K^T-Y=GG^T
    \end{equation}
\end{mylem}
\quad \textit{Proof}: Consider the increment of \eqref{X}
\begin{equation*}
    \begin{split}
        -dK^TB^TXA_K + A_K^TdXA_K - A_K^TXBdK-dX \\ 
        -dK^TD^TC_K - C_K^TDdK=0
    \end{split}
\end{equation*}
Pre-multiplying both sides with $Y$ and take trace, we get
\begin{equation*}
\begin{split}
    &dJ(K) = Tr(GG^TdX) = \\
    &2Tr(B^TXA_KYdK^T)+2Tr(D^TC_KYdK^T)
\end{split}
\end{equation*}
Reorganizing the terms, we obtain
\begin{equation*}
    dJ(K)=2Tr((B^TXA_K+D^TC_K)YdK^T)
\end{equation*}
and thus $\nabla J(K) = 2(B^TXA_K-D^TDK)Y$ since $D^TC=0$. The proof is complete. $\hfill \square$

In this paper, we will solve the Problem DHSCCS from an optimization perspective. The Lagrangian function of the problem in \eqref{equivalent problem} is constructed as follows.
\begin{equation*}
    \label{lagrangian}
    L(K,Q) = J\left( K \right) -Tr\left( Q^T\left[ \begin{array}{c}
	A-BK\\
	\mathbf{1}^T-\mathbf{1}^T\left( A-BK \right)\\
\end{array} \right] \right) 
\end{equation*}
where $Q$ is known as the Lagrangian multiplier. Due to the non-convex nature of the problem, a global optimal solution is often unattainable. Consequently, we will focus on local optimality instead of global optimality. In what follows, we illustrate how the necessary conditions for local optimality can be deduced using the KKT conditions.  Considering the problem outlined in \eqref{equivalent problem}, the inequality constraints are affine with respect to $K$, and thus it satisfies the linearity constraint qualification (LCQ) \cite{eustaquio2009karush}. This compliance leads to strong duality, ensuring a duality gap of $0$. In other words, if $K^*$ is the local optimal solution to \eqref{equivalent problem}, there must exist a dual variable $Q^*$, known as KKT multiplier, such that the following KKT conditions hold.
%\begin{equation}
%    \label{KKT conditions}
%    \begin{cases}
\begin{equation}
\label{KKT conditions}
    \begin{aligned}
        &\left\{ \begin{array}{l}
	\nabla J\left( K^* \right) -\left[ \begin{matrix}
	-B^T&		B^T \mathbf{1}\\
\end{matrix} \right] Q^*=0\\
	Q^*\ge 0\\
	\left[ \begin{array}{c}
	A-BK^*\\
	\mathbf{1}^T-\mathbf{1}^T\left( A-BK^* \right)\\
\end{array} \right] \ge 0\\
	Tr\left( Q^{*^T}\left[ \begin{array}{c}
	A-BK^*\\
	\mathbf{1}^T-\mathbf{1}^T\left( A-BK^* \right)\\
\end{array} \right] \right)=0\\
        \end{array} 
        \right.
    \end{aligned}
\end{equation}
%\end{cases}
%\end{equation}
\begin{mydef}
$K\in \mathcal{F}$ is called a stationary point of problem \eqref{equivalent problem} if it satisfies its KKT conditions in \eqref{KKT conditions}.  
\end{mydef}
\begin{myrem}
    Notice that the local optimal solution to \eqref{equivalent problem} must be a stationary point. However, a stationary point is not necessarily a locally optimal point due to the existence of saddle points. Although getting sufficient conditions for local optimality is challenging, we can still eliminate points that are not locally optimal or saddle. The stationary point is usually the highest pursuit for non-convex optimization problems. 
\end{myrem}
\section{First-Order Method}\label{First-Order Method}
To the best of our knowledge, in $H_2$ control setting, most existing works deal with inequality constraints via dual methods, including the alternating direction method of multipliers (ADMM) and the augmented Lagrangian Method (ALM). However, although they have shown great performance, the convergence of these methods is not guaranteed in most non-convex cases. Other first-order methods such as projected
gradient descent (PGD) could lead to infeasible projection points, since the boundary points of \eqref{equivalent problem} may not be a stabilizing gain \cite{valcher2018state}. In this section, we propose a first-order interior point method (FIPM) to solve \eqref{equivalent problem}. We will show that FIPM converges to a stationary point of the original problem. We first reformulate \eqref{equivalent problem} into the following unconstrained problem.
\begin{equation}
\label{FIPM first order problem}
    \min_K J(K)-\cfrac{1}{t}\sum_{i,j}\mathrm{log}\left[ \begin{array}{c}
	A-BK\\
	\mathbf{1}^T-\mathbf{1}^T\left( A-BK \right)\\
    \end{array} \right]_{ij}
\end{equation}
where $t>0$ is the penalty parameter. The intuition for using the log-barrier function is that, as $t$ approaches $+\infty$, the log-barrier function will approximate an indicator function. Therefore, the solution $K$, as $t$ approaches $+\infty$, will converge to the solution of the original problem. Denote
\begin{equation}
\label{JLBF}
\begin{split}
        &JLBF(K,t) = J(K)\\
        &-\cfrac{1}{t}\sum_{i,j}\mathrm{log}\left(\left[ \begin{array}{c}
	A\\
	\mathbf{1}^T-\mathbf{1}^TA\\
\end{array} \right] +\left[ \begin{array}{c}
	-B\\
	\mathbf{1}^TB\\
\end{array} \right] K\right)_{ij}
\end{split}
\end{equation}
Since \eqref{FIPM first order problem} is an unconstrained problem and the object function is smooth, we use gradient descent to solve $\eqref{FIPM first order problem}$. The gradient of $JLBF$ is provided below. The full algorithm is shown in \cref{Alg: 1}.
\begin{equation}
\label{nabla JLBF}
\begin{split}
        &\nabla_K JLBF(K,t) = \nabla J(K)\\
        &-\cfrac{1}{t}\cdot\left[\begin{matrix}
	-B^T&		B^T \mathbf{1}\\
\end{matrix} \right]\cdot \left(\left[\begin{array}{c}A_K\\\mathbf{1}^T-\mathbf{1}^TA_K\\\end{array} \right]\right)^{\circ-1}
\end{split}
\end{equation}

\begin{mypro}
\label{convergence of algorithm 1}
    \cref{Alg: 1} converges to a stationary point of problem \eqref{equivalent problem}.
\end{mypro}

\begin{algorithm}[H]
\SetAlgoLined
\textbf{Input:} $K^{(0)}$, $t^{(0)}, \mu>1$\\
\textbf{Output:} $K^*$\\
Initialize $h=0$;\\
\Repeat{ $||K^{(h+1)}-K^{(h)}||<\epsilon_2$\vspace{5pt}}{
\vspace{0pt}
\vspace{0pt}
    - Initialize $k=0$; \\
    - $\hat{K}^{(k)} = K^{(h)}$;\\
    \Repeat{$||\nabla_K JLBF(\hat{K}^{(k)},t^{(h)})||<\epsilon_1$ or $\hat{K}^{(k)}=\hat{K}^{(k-1)}$\vspace{5pt}}{
    \vspace{2pt}
    - Calculate $\nabla_K JLBF(\hat{K}^{(k)},t^{(h)})$ by \eqref{nabla JLBF};\\
    - $\hat{K}^{(k+1)} = \hat{K}^{(k)} - s\nabla_K JLBF(\hat{K}^{(k)},t^{(h)})$, where $s$ is determined by Armijo rule \cite{nocedal1999numerical};\\
    - $k = k + 1$;\\
    \vspace{-5pt}
    }
    - $K^{(h+1)} = \hat{K}^{(k)}$;\\
    - $t^{(h+1)}=\mu t^{(h)}$;\\
    - $h = h + 1$;\\
    \vspace{-5pt}}
\Return $K^{(h)}$.
\caption{FIPM for solving \textbf{Problem DHSCCS}}\label{Alg: 1}
\end{algorithm}
\quad \textit{Proof}: Given the continuity and coerciveness of $JLBF$, according to Weierstrass' Theorem, $JLBF$ attains the global minimum within its domain, establishing its lower-boundedness. Employing gradient descent, we generate a non-increasing sequence, ensuring the convergence of $\hat{K}^{(k)}$ for each penalty parameter $t^{(h)}$. Additionally, as $t\rightarrow \infty$, the FIPM solution converges to the solution of the original problem \eqref{equivalent problem}. The proof of convergence is complete. Now we move on to the stationary part. Since the algorithm is guaranteed to converge, we obtain 
\begin{equation}
\begin{split}
        &\lim_{h\rightarrow \infty}\nabla_K JLBF(K^{(h+1)},t^{(h)}) = \nabla J(K^*)\\
        &-\left[\begin{matrix}
	-B^T&		B^T \mathbf{1}\\
\end{matrix} \right]\cdot Q^*
\end{split}
\end{equation}
where $Q^* = \lim_{h\rightarrow \infty}\cfrac{1}{t^{(h)}}\cdot \left(\left[\begin{array}{c}A_K^{(h+1)}\\\mathbf{1}^T-\mathbf{1}^TA_K^{(h+1)}\\\end{array} \right]\right)^{\circ-1}\ge0$ and $A_K^{(h+1)}=A-BK^{(h+1)}$ because $K$ always exist in the feasible region and $t^{(h)}>0$. Hence, the first three conditions of \eqref{KKT conditions} hold trivially. For the last condition, denote $S^*=\lim_{h\rightarrow \infty}\left[\begin{array}{c}A_K^{(h+1)}\\\mathbf{1}^T-\mathbf{1}^TA_K^{(h+1)}\\\end{array} \right]$. If there is an element $S^*_{ij}>0$, $Q^*_{ij}=0$ since $1/t^{(h)}$ will force it to $0$ as $h\rightarrow \infty$. Since the condition holds trivially, we do not need to consider $S^*_{ij}=0$. The proof is complete. $\hfill \square$

\section{Second-Order Method}\label{Second-Order Method}
In this section, we propose a second-order interior 
point method (SIPM) to solve \eqref{equivalent problem}. 
The main challenges involve calculating the Hessian matrix for $J(K)$ and the second log-barrier term in \eqref{JLBF}. 
Since both terms are scalar functions dependent on a matrix, their Hessian matrices exist in a four-dimension space. To avoid discussing tensors, we vectorize $K$, ensuring that the Hessian matrices are represented in two dimensions for ease of mathematical expression.
To facilitate discussions, denote 
\begin{equation*}
    \begin{aligned}
        &\varGamma X=-C_K^TC_K\\
        &\varGamma^*Y=GG^T\\
        &\varGamma_{mn}X = \partial_{K_{mn}} A_K^TXA_K+A_K^TX\partial_{K_{mn}} A_K\\
        &\varGamma_{mn}^* Y=\partial_{K_{mn}} A_K Y A_K^T + A_K Y \partial_{K_{mn}} A_K^T
    \end{aligned}
\end{equation*}
where we define the general Lyapunov operator $\varGamma:\mathbb{R}^{n\times n}\rightarrow \mathbb{R}^{n\times n}$ as $P\mapsto A_K^T P A_K - P$, and $K_{mn}$ denotes the entry on $m$-th row and $n$-th column of $K$. Furthermore, we have the inner product property $\left< \varGamma X, Y \right>=\left< X, \varGamma^*Y \right>  $.
\subsection{Derivation of Hessian matrix}
In this subsection, we explicitly compute the Hessian matrix of $J(K)$. We first express the Hessian matrix of $J(K)$ as the partial derivative of $\nabla J(K)$ for each entry of $K$, and then we obtain $\nabla_{\mathrm{vec}(K)}^2 J(K)$. Similarly, we derive the gradient and Hessian matrix of the log barrier term in a vectorized sense. 
% With the gradient and Hessian matrix of $JLBF$, we are ready to propose the SIPM algorithm.
\begin{mythr}
Consider the system \eqref{rewritten system}, the partial derivative of $\nabla J(K)$ to each entry of $K$ can be obtained as
\begin{equation}
\label{Hessian J(K)}
\begin{split}
    &\partial_K\partial_{K_{ij}}J = -2B^T\left[X^{ij}+\left(X^{ij}\right)^T\right]\left(A-BK\right)Y\\
    &-2\left[B^TX\left(A-BK\right)-D^TDK\right]\left[Y^{ij}+\left(Y^{ij}\right)^T\right]\\
    &-2\left[B^TXB+D^TD\right]S^{ij}Y\\
    &+2B^T\left[Z^{ij}+\left(Z^{ij}\right)^T\right]\left(A-BK\right)Y
\end{split}
\end{equation}
where $S^{ij}$ denotes a sparse matrix with $1$ on $i$-th row and $j$-th column and $0$ on all other entries.
\begin{subequations}
    \begin{align*}
        & \varGamma X^{ij} = -\left(A-BK\right)^TXBS^{ij} \\
        & \varGamma^*Y^{ij} = -\left(A-BK\right)Y\left(S^{ij}\right)^TB^T \\
        & \varGamma Z^{ij} = -K^TD^TDS^{ij}
    \end{align*}
\end{subequations}
\end{mythr}
\quad \textit{Proof}: From \eqref{nabla J(K)}, the gradient of $J$ with a single entry of $K$ can be expressed as 
\begin{equation*}
\begin{split}
    \partial_{K_{ij}}J &=2\left<B^TXA_KY,S^{ij}\right>-2\left<D^TDKY,S^{ij}\right>\\
    &=2\left<Y,A_K^TXBS^{ij}\right>-2\left<Y,K^TD^TDS^{ij}\right>
\end{split}
\end{equation*}
Therefore, each entry of the Hessian matrix is 
\begin{equation*}
\begin{split}
    &\partial_{K_{mn}}\partial_{K_{ij}}J=2\left<\partial_{K_{mn}}Y,A_K^TXBS^{ij}\right>\\
    &-2\left<Y,\left(BS^{mn}\right)^TXBS^{ij}\right>+2\left<Y,A_K^T\partial_{K_{mn}}XBS^{ij}\right>\\
    &-2\left<\partial_{K_{mn}}Y,K^TD^TDS^{ij}\right>-2\left<Y,\left(S^{mn}\right)^TD^TDS^{ij}\right>
\end{split}
\end{equation*}
By using the inner product property of $\varGamma$ and $\varGamma^*$, we have
\begin{equation*}
\begin{split}
    &\partial_{K_{mn}}\partial_{K_{ij}}J=2\left<\varGamma^*_{mn}Y,\varGamma^{-1}\left[-\left(A-BK\right)^TXBS^{ij}\right]\right>\\
    &-2\left<Y,A_K^T\left[\varGamma^{-1}\left(\varGamma_{mn}X\right)+\varGamma^{-1}\left(\partial_{K_{mn}}C_K^TC_K\right)\right]BS^{ij}\right>\\
    &+2\left<Y,\left(-BS^{mn}H\right)^TXBS^{ij}\right>-2\left<Y,\left(S^{mn}\right)^TD^TDS^{ij}\right>\\
    &-2\left<\varGamma^*_{mn}Y,\varGamma^{-1}\left(-K^TD^TDS^{ij}\right)\right>\\
\end{split}
\end{equation*}
From $\varGamma_{mn}X = \left(-BS^{mn}\right)^TXA_K+A_K^TX\left(-BS^{mn}\right)$ and $\varGamma^*_{mn}Y=\left(-BS^{mn}\right)YA_K^T+A_KY\left(-BS^{mn}\right)^T$, we have 
\begin{equation*}
\begin{split}
    &\partial_{K_{mn}}\partial_{K_{ij}}J\\
    &=2\left<\left(-BS^{mn}\right)YA_K^T+A_KY\left(-BS^{mn}\right)^T,X^{ij}\right>\\
    &+2\left<Y^{ij},\left(-BS^{mn}\right)^TXA_K+A_K^TX\left(-BS^{mn}\right)\right>\\
    &+\left<Y^{ij},\left(-DS^{mn}\right)^TC_K+C_K^T\left(-DS^{mn}\right)\right>\\
    &+2\left<Y,\left(-BS^{mn}\right)^TXBS^{ij}\right>-2\left<Y,\left(S^{mn}\right)^TD^TDS^{ij}\right>\\
    &-2\left<\left(-BS^{mn}\right)YA_K^T+A_KY\left(-BS^{mn}\right)^T,Z^{ij}\right>\\
\end{split}
\end{equation*}
By moving $S^{mn}$ and $(S^{mn})^T$ terms to one side, we have
\begin{equation*}
\begin{split}
    &\partial_{K_{mn}}\partial_{K_{ij}}J=\\
    &2\left<S^{mn},-B^TX^{ij}A_KY\right>-2\left<\left(S^{mn}\right)^T,YA_K^TX^{ij}B\right>\\
    &-2\left<Y^{ij}A_K^TXB,\left(S^{mn}\right)^T\right>-2\left<B^TXA_KY^{ij},S^{mn}\right>\\
    &+2\left<-Y^{ij}C_K^TD,\left(S^{mn}\right)^T\right>-2\left<D^TC_KY^{ij},S^{mn}\right>\\
    &+2\left<-Y\left(S^{ij}\right)^TB^TXB,\left(S^{mn}\right)^T\right>\\
    &+2\left<S^{mn},B^TZ^{ij}A_KY\right>+2\left<\left(S^{mn}\right)^T,YA_K^TZ^{ij}B\right>\\
    &-2\left<Y\left(S^{ij}\right)^TD^TD,\left(S^{mn}\right)^T\right>
\end{split}
\end{equation*}
Now we can extend $\partial_{K_{mn}}$ to $\partial_K$ by remove the $S^{mn}$ and $(S^{mn})^T$ terms. After reorganizing terms and combining similar terms we obtain \eqref{Hessian J(K)}. The proof is complete.$\hfill \square$

If we regard $K$ as $\mathrm{vec}(K)$ instead of a matrix, by leveraging \eqref{Hessian J(K)}, we can obtain the gradient and Hessian expression of $J(K)$ as
\begin{equation}
\label{G_J}
    G_J(K) = \mathrm{vec}(\nabla J(K))
\end{equation}
\begin{equation}
\label{H_J}
\begin{split}
    H_J(K) = &\left[\mathrm{vec}\left(\cfrac{\partial^2J}{\partial_{K_{11}}\partial_K}\right) ,\ldots, \mathrm{vec}\left(\cfrac{\partial^2J}{\partial_{K_{m1}}\partial_K}\right),\ldots\right.\\ 
    &\left. \mathrm{vec}\left(\cfrac{\partial^2J}{\partial_{K_{1p}}\partial_K}\right) ,\ldots, \mathrm{vec}\left(\cfrac{\partial^2J}{\partial_{K_{mp}}\partial_K}\right)\right]
\end{split}
\end{equation}
After addressing the $J(K)$, we proceed to discuss the log-barrier term. Denote
\begin{equation}
    LBF(K,t) = -\cfrac{1}{t}\sum_{i,j}\mathrm{log}\left[ \begin{array}{c}
	A-BK\\
	\mathbf{1}^T-\mathbf{1}^T\left( A-BK \right)\\
    \end{array} \right]_{ij}
\end{equation}
Again we regard $K$ as $\mathrm{vec}(K)$ instead of a matrix. We propose \cref{nabla hessian LBF} to express the gradient and Hessian matrix of $LBF(K,t)$. Before that, we provide a technical lemma which is critical in the proof of \cref{nabla hessian LBF}.
\begin{mylem}
\label{technical equality}
    Given three general matrices $\mathcal{A},\mathcal{X},\mathcal{B}$ with compatible dimensions, the following equality holds.
    \begin{equation}
        \mathrm{vec}(\mathcal{A}\mathcal{X}\mathcal{B}) = (\mathcal{B}^T\otimes \mathcal{A})\cdot \mathrm{vec}(\mathcal{X})
    \end{equation}
\end{mylem}
% \quad \textit{Proof}: We can first express $\mathcal{AX}$ by columns, and each column of $\mathcal{AXB}$ can be obtained as the product of $\mathcal{AX}$ with each column of $\mathcal{B}$. After vectorization, we can extract $\mathrm{vec}(\mathcal{X})$ and the coefficient matrix happens to be $\mathcal{B}^T\otimes A$.$\hfill \square$
    \begin{algorithm}[H]
\SetAlgoLined
\textbf{Input:} $K^{(0)}$, $t^{(0)}, \mu>1$\\
\textbf{Output:} $K^*$\\
Initialize $h=0$;\\
\Repeat{ $||K^{(h+1)}-K^{(h)}||<\epsilon_2$\vspace{5pt}}{
\vspace{0pt}
\vspace{0pt}
    - Initialize $k=0$; \\
    - $\hat{K}^{(k)} = K^{(h)}$;\\
    \Repeat{$||G_{JLBF}(\hat{K}^{(k)},t^{(h)})||<\epsilon_1$ or $\hat{K}^{(k)}=\hat{K}^{(k-1)}$\vspace{5pt}}{
    \vspace{2pt}
    - Calculate the gradient $G_{JLBF}(\hat{K}^{(k)},t^{(h)})$ from \eqref{G_J}, \eqref{G LBF} and \eqref{G};\\
    - Calculate the Hessian $H_{JLBF}(\hat{K}^{(k)},t^{(h)})$ from \eqref{H_J}, \eqref{H LBF} and \eqref{H};\\
    - Calculate the modified Hessian from \eqref{modified Hessian};\\
    - Calculate the modified Newton step from \eqref{modified newton step};\\
    - $\hat{K}^{(k+1)} = \hat{K}^{(k)} - s\cdot \mathrm{mat}(\mathrm{vec}(\Delta K))$, where we first recover the modified Newton step to a matrix. The $s$ is determined by Armijo rule \cite{nocedal1999numerical} with initial value 1;\\
    - $k = k + 1$;\\
    \vspace{-5pt}
    }
    - $K^{(h+1)} = \hat{K}^{(k)}$;\\
    - $t^{(h+1)}=\mu t^{(h)}$;\\
    - $h = h + 1$;\\
    \vspace{-5pt}}
\Return $K^{(h)}$.
\caption{SIPM for solving \textbf{Problem DHSCCS}}\label{Alg: 2}
\end{algorithm}
\begin{mythr}
    The gradient and Hessian matrix of $LBF$ in a $\mathrm{vec}(K)$ sense can be expressed as follows
    \label{nabla hessian LBF}
    \begin{equation}
    \label{G LBF}
    \begin{split}
        &G_{LBF}(K,t)= -\cfrac{1}{t}\cdot I \otimes \left[\begin{matrix}
	-B^T&		B^T \mathbf{1}\\
\end{matrix} \right] \\
&\cdot \left\{ \mathrm{vec}\left(\left[ \begin{array}{c}
	A\\
	\mathbf{1}^T-\mathbf{1}^TA\\
\end{array} \right]\right)       \right.   \left. +I^T\otimes \left[ \begin{array}{c}
	-B\\
	\mathbf{1}^TB\\
\end{array} \right]\cdot \mathrm{vec}(K)\right\}^{\circ-1}
    \end{split}
    \end{equation}
    \begin{equation}
    \label{H LBF}
        \begin{split}
        &H_{LBF}(K,t) = \cfrac{1}{t}\cdot I \otimes \left[\begin{matrix}
	-B^T&		B^T \mathbf{1}\\
\end{matrix} \right] \\
        &\cdot \mathrm{diag}\bigg(  \bigg\{  \mathrm{vec}\left(\left[ \begin{array}{c}
	A\\
	\mathbf{1}^T-\mathbf{1}^TA\\
\end{array} \right]\right)  +   I^T\otimes \left[ \begin{array}{c}
	-B\\
	\mathbf{1}^TB\\
\end{array} \right]  \\
    & \cdot \mathrm{vec}(K) \bigg\}^{\circ-2}     \bigg) \cdot I^T \otimes \left[ \begin{array}{c}
	-B\\
	\mathbf{1}^TB\\
\end{array} \right]
        \end{split}
    \end{equation}
\end{mythr}
\quad \textit{Proof}: For the simplicity, we denote $\mathcal{A}=\left[ \begin{array}{c}
	A\\
	\mathbf{1}^T-\mathbf{1}^TA\\
\end{array} \right]$, $\mathcal{B}=\left[ \begin{array}{c}
	-B\\
	\mathbf{1}^TB\\
\end{array} \right]$, $\mathcal{C} = I$. Then 
\begin{equation*}
    \nabla_K\sum_{ij}\mathrm{log}(\mathcal{A}+\mathcal{B}K\mathcal{C})_{ij}=\mathcal{B}^T(\mathcal{A}+\mathcal{B}K\mathcal{C})^{\circ-1}\mathcal{C}^T
\end{equation*}
Use \cref{technical equality} twice,
\begin{equation*}
    \begin{split}
    &\mathrm{vec}\left(\mathcal{B}^T(\mathcal{A}+\mathcal{B}K\mathcal{C})^{\circ-1}\mathcal{C}^T\right)\\
    &=\left(I\otimes \mathcal{B}^T\right)\cdot \left[\mathrm{vec}(\mathcal{A})+\mathrm{vec}(\mathcal{B}K\mathcal{C})\right]^{\circ-1} \\
    &=\left(I\otimes \mathcal{B}^T\right)\cdot \left[\mathrm{vec}(\mathcal{A}) + (I^T\otimes \mathcal{B})\cdot \mathrm{vec}(K)\right]^{\circ-1}
    \end{split}
\end{equation*}
By recovering $\mathcal{A,B}$, we get \eqref{G LBF}. Then \eqref{H LBF} is trivial by using chain rules. The proof is complete.$\hfill \square$
\subsection{Hessian Modification}
The indefiniteness of the Hessian matrix in non-convex problems presents a common challenge, where the Newton step may not ensure a descent direction. Even when the Hessian matrix is positive definite and the Newton step is a descent direction, the algorithm could still converge to a saddle point. Presently, a universal solution to this issue is still unresolved. In other words, evading saddle points in non-convex problems is still an open problem.
% Some existing works have proposed a few heuristic approaches to guarantee the positive definiteness of the Hessian matrix. Paternain \cite{paternain2019newton} proposed to replace all negative eigenvalues with their absolute values. 
% Nocedal \cite{nocedal1999numerical} made the Hessian matrix sufficiently positive definite by adding a small matrix. 
Thus, our goal is to force the Hessian matrix to be positive definite to ensure that the Newton step represents a descent direction.
In this paper, we adopt the diagonal modification method proposed by Nocedal \cite{nocedal1999numerical}. This method sets a lower bound for the eigenvalues, thus ensuring the Hessian matrix's positive definiteness. Consequently, the modified Newton step is a descent direction and secure convergence towards a sub-optimal solution.
\subsection{SIPM}
In this subsection, we will provide the full algorithm of SIPM. The vectorized gradient and the Hessian matrix of $JLBF$ are provided as follows
\begin{subequations}
\begin{align}
    &G_{JLBF}(K,t) = G_J(K) + G_{LBF}(K,t) \label{G}\\
    &H_{JLBF}(K,t) = H_J(K) + H_{LBF}(K,t) \label{H}
\end{align}
\end{subequations}
% Given $t$, the second-order Taylor approximation of $JLBF$ around $K$ is
% \begin{equation}
% \begin{split}
%     &\widehat{JLBF}(K+\Delta K, t) \\
%     & = JLBF(K,t) + G_{JLBF}(K,t)^T\cdot \mathrm{vec}(\Delta K) \\
%     &+ \cfrac{1}{2}\cdot \mathrm{vec}(\Delta K)^T \cdot H_{JLBF}(K,t) \cdot \mathrm{vec}(\Delta K)
% \end{split}
% \end{equation}
Before deriving the Newton step, we will first implement the diagonal modification method. The spectral decomposition of $H_{JLBF}$ is $H_{JLBF}=Q\Sigma Q^T$, where $Q$ is an orthogonal matrix and $\Sigma$ is a diagonal matrix with the eigenvalues on its diagonal. The modified Hessian matrix can be expressed as 
\begin{equation}
\label{modified Hessian}
H_{JLBF_\delta}=Q\Sigma_\delta Q^T
\end{equation}
where
\begin{equation}
    (\Sigma_\delta)_{ii}= \begin{cases}(\Sigma)_{ii}, &(\Sigma)_{ii} \geq \delta \\ \delta, &(\Sigma)_{ii}<\delta\end{cases}
\end{equation}
The modified Newton step can now be expressed as
\begin{equation}
\label{modified newton step}
    \mathrm{vec}(\Delta K)=H_{JLBF_\delta}(K,t)^{-1}\cdot G_{JLBF}(K,t)
\end{equation}
The full algorithm is shown in \cref{Alg: 2}. 
% \begin{algorithm}[H]
% \SetAlgoLined
% \textbf{Input:} $K^{(0)}$, $t^{(0)}, \mu>1$\\
% \textbf{Output:} $K^*$\\
% Initialize $h=0$;\\
% \Repeat{ $||K^{(h+1)}-K^{(h)}||<\epsilon_2$\vspace{5pt}}{
% \vspace{0pt}
% \vspace{0pt}
%     - Initialize $k=0$; \\
%     - $\hat{K}^{(k)} = K^{(h)}$;\\
%     \Repeat{$||G_{JLBF}(\hat{K}^{(k)},t^{(h)})||<\epsilon_1$ or $\hat{K}^{(k)}=\hat{K}^{(k-1)}$\vspace{5pt}}{
%     \vspace{2pt}
%     - Calculate the gradient $G_{JLBF}(\hat{K}^{(k)},t^{(h)})$ from \eqref{G_J}, \eqref{G LBF} and \eqref{G};\\
%     - Calculate the Hessian $H_{JLBF}(\hat{K}^{(k)},t^{(h)})$ from \eqref{H_J}, \eqref{H LBF} and \eqref{H};\\
%     - Calculate the modified Hessian from \eqref{modified Hessian};\\
%     - Calculate the modified Newton step from \eqref{modified newton step};\\
%     - $\hat{K}^{(k+1)} = \hat{K}^{(k)} - s\cdot \mathrm{mat}(\mathrm{vec}(\Delta K))$, where we first recover the modified Newton step to a matrix. The $s$ is determined by Armijo rule\cite{nocedal1999numerical} with initial value 1;\\
%     - $k = k + 1$;\\
%     \vspace{-5pt}
%     }
%     - $K^{(h+1)} = \hat{K}^{(k)}$;\\
%     - $t^{(h+1)}=\mu t^{(h)}$;\\
%     - $h = h + 1$;\\
%     \vspace{-5pt}}
% \Return $K^{(h)}$.
% \caption{SIPM for solving \textbf{Problem DHSCCS}}\label{Alg: 2}
% \end{algorithm}

\begin{mypro}
    \cref{Alg: 2} converges to a stationary point of problem \eqref{equivalent problem}.
\end{mypro}

\quad \textit{Proof}: It is similar to \cref{convergence of algorithm 1} and is omitted.$\hfill \square$
\begin{myrem}
    Throughout the previous discussions, we assume we can initialize a strictly feasible controller $K^{(0)}$, i.e., $\left[\begin{array}{c}A_{K^{(0)}}\\\mathbf{1}^T-\mathbf{1}^TA_{K^{(0)}}\\\end{array} \right]>0$ element-wisely, since otherwise the $\mathrm{log}(\cdot)$ will be infinite and the derivative does not exist. However, this inequality might not be strict for a general $K^{(0)}$ as the equality could hold for some stabilizing gains \cite{valcher2018state}. A trick that avoids this issue is to relax the constraint by a sufficiently small number $\epsilon_r>0$. In other words, the constraint of the original problem \eqref{equivalent problem} should be modified to $\left[\begin{array}{c}A_K\\\mathbf{1}^T-\mathbf{1}^TA_K\\\end{array} \right]\ge -\epsilon_r$. The initialization of $K^{(0)}$ is discussed in the next section. Our methods can easily be modified to consider this numerical modification. The $\epsilon_r$ is omitted in our main analysis to avoid confusion.
\end{myrem}

\begin{myrem}
    In terms of the complexity of FIPM and SIPM. Since 
    $J(K)$ is strongly convex around the global minimizer \cite{bu2019lqr} and the log barrier term is convex, we can estimate the iteration complexity of FIPM and SIPM by assuming they are 
    linearly convergent and quadratically convergent respectively.
    % , of which the iteration complexities are given by $\mathcal{O}(\text{log}(1/\epsilon_1^2))$ and $\mathcal{O}(\text{log}(\text{log}(1/\epsilon_1^2)))$. 
    During each iteration, the computational complexity of deriving the gradient and Hessian are $\mathcal{O}(n^3)$ and $\mathcal{O}(mn^4)$. The computation of the modified Newton step involves the inverse of the Hessian matrix, which costs $\mathcal{O}(m^3n^3)$. 
    The speed of FIPM and SIPM balances between iteration numbers and computational costs per iteration.
\end{myrem}

\section{Discussions on the initialization of $K^{(0)}$}\label{Discussions on the initialization of K0}
The interior point method requires the $K^{(0)}$ to be strictly feasible. More specifically, we need to find $K^{(0)}$ to make $A-BK^{(0)}$ compartmental and Schur. We first introduce a useful lemma to facilitate discussions.
\begin{mylem}
\label{schur and hurwitz}
    \cite{liu2023proportional} For a non-negative matrix $A-BK$, the following statements are equivalent:
    
    1) Matrix $A-BK$ is Schur;

    2) Matrix $A-BK-I$ is Hurwitz;
\end{mylem}
\cref{schur and hurwitz} associates Schur and Hurwitz, allowing us to design a controller that makes $A-BK$ Schur and compartmental by indirectly finding such a controller that makes $A-I-BK$ Hurwitz with compartmental $A-BK$. Here we borrow the idea, concepts and notations from \cite{valcher2018state}. The procedure to find such $K^{(0)}$ is listed below
\begin{itemize}
    \item If $A-I$ is irreducible, we use Lemma 7 in \cite{valcher2018state}. If there is a solution $K$ with sufficiently small $\epsilon>0$ s.t. $A-BK$ is compartmental, we can use this $K$ as $K^{(0)}$.
    \item If $A-I$ is reducible, we use Proposition 15 in \cite{valcher2018state}. If there is a solution $K$ with sufficiently small $\epsilon>0$ s.t. $A-BK$ is compartmental, we can use this $K$ as $K^{(0)}$.
\end{itemize}

\section{Simulations}\label{Simulations}

\begin{table*}[t]
\caption{Comparison of FIPM, SIPM and ADMM with increasing system scales for example 1}
\centering

\begin{tabular}{|l|r|r|r|r|r|r|r|r|r|r|}
\hline
\label{increasing dimension example 1}
 \diagbox{Method}{$N$}& 1 & 2 & 3 & 4 & 5 & 6 & 7 & 8 & 9 & 10 \\ \hline 
 FIPM& 7.2385s & 16.2040s & 23.3113s & 40.2485s & 46.7011s & 57.3487s & 71.2840s & 84.9379s & 96.7251s & 102.6486s \\ \hline
 SIPM& 0.0862s & 0.5602s & 1.6532s & 4.9985s & 9.4232s & 32.6799s & 57.8460s & 92.0200s & 165.9687s & 225.2351s \\ \hline
 ADMM \cite{yang2024h2}& 0.4511s & 2.4909s & 5.1267s & 10.1377s & 23.6348s & 38.7567s & 69.9793s & 77.3889s & 129.4604s & 176.6897s \\ \hline
\end{tabular}

\end{table*}
\begin{figure}[t]
    \centering    \includegraphics[width=0.5\linewidth]{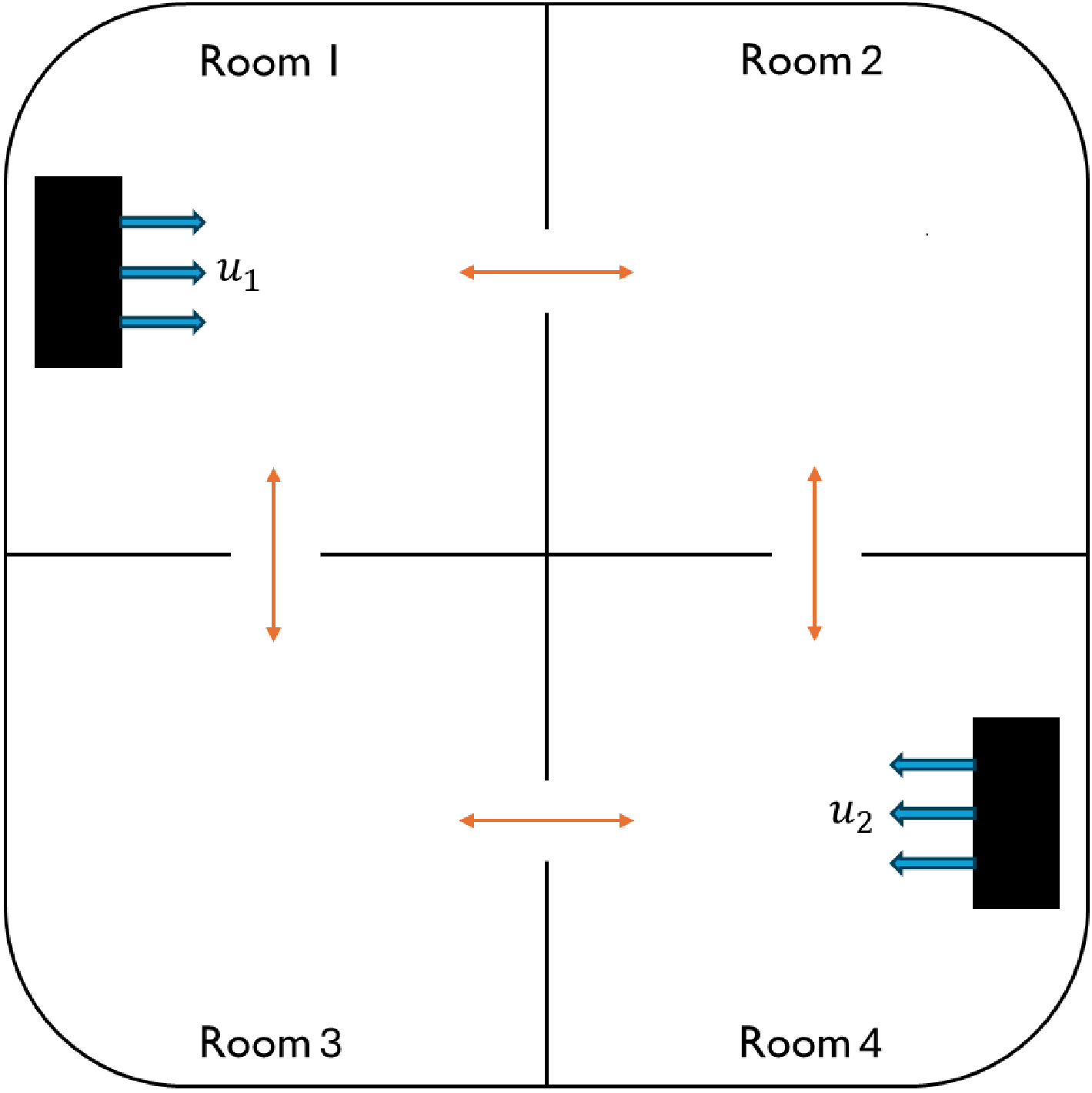}
    \caption{A 4-room thermal system with two inputs}
    \label{room temperature}
\end{figure}
% In this section, we will provide some numerical simulations of a real-world thermal system to compare the FIPM and SIPM concerning, 1) the performance in each external iteration when there is one thermal system, and 2) the overall performance with the increasing number of thermal systems. 
\subsection{Example 1: Thermal System}
\cref{room temperature} illustrates the thermal system, with no heat exchange with the environment and two rooms directly connected to the heat conditioner. The model can be represented as $\{\mathcal{V},\mathcal{E}\}$, where $\mathcal{V}$ represents the set of rooms, and $\mathcal{E}$ represents the set of directed heat flow. The red double-sided arrow represents the heat exchange between the adjacent rooms. The heat of this thermal system is depicted by an ordinary differential equation (ODE) with the heat in room $i$ given by 
\begin{equation*}
    \begin{split}
        \Dot{x}_i(t) = a_ix_i(t)+\sum_{ (j,i) \in \mathcal{E}}{a_{ij}x_j(t)}-\sum_{(i,j) \in \mathcal{E}}{a_{ji}x_i(t)}+\\
        \sum_{k}{b_{ik}u_k(t)}+w_i(t)
    \end{split}
\end{equation*}
where $a_ix_i(t)$ represents the decay, $a_{ij}x_j(t)$ represents the heat flow from room $j$ to room $i$, $a_{ji}x_i(t)$ represents the heat flow from room $i$ to room $j$, $b_{ik}u_k(t)$ represents the input from actuator $k$ to room $i$, $w_i(t)$ is the disturbance. We assume that we can not measure heat continuously due to the digitalization of devices. Therefore, we sample the measurement with interval $\Delta t=0.1s$ and the system can be rewritten as \eqref{basic system} with system matrices listed as follows
\begin{equation*}
    \begin{split}
        &A = \left[ \begin{matrix}
	0.5&		0.2&		0.1&		0.0\\
	0.1&		0.6&		0.0&		0.2\\
	0.4&		0.0&		0.8&		0.4\\
	0.0&		0.2&		0.1&		0.4\\
\end{matrix} \right], B = \left[ \begin{matrix}
	0.1&		0.0\\
	0.0&		0.0\\
	0.0&		0.0\\
	0.0&		0.1\\
\end{matrix} \right], G = I_{4\times4},\\
&C = \left[ \begin{matrix}
	I_{2\times2}&		I_{2\times2}\\
	0_{2\times2}&		0_{2\times2}\\
\end{matrix} \right], D = \left[ \begin{matrix}
	0_{2\times2}\\
	I_{2\times2}\\
\end{matrix} \right]
    \end{split}
\end{equation*}
with $C^TD=0$, $D^TD\succ0$. Physically, our goal is to design a state-feedback controller that transforms an unstable system, which operates without heat leakage, into a stable system that accounts for heat leakage, while simultaneously enhancing the thermal system's robustness. Mathematically, we aim to design a state-feedback controller $K$ to minimize the system's $H_2$ norm, ensuring that $A_K$ is both compartmental and Schur. We run FIPM and SIPM for 10 external iterations separately to compare their performance. 

We choose $\epsilon_1=1e-4$, $\epsilon_2=1e-3$, $\epsilon_r=1e-9$, $\delta=1$, $\mu = 4$ for FIPM and SIPM. We can easily check $A-I$ is irreducible since $A$ describes a strongly-connected graph. By leveraging Lemma 7 in \cite{valcher2018state} and take $\mathbf{v}_1=\left[ \begin{matrix}
	-4& 1
\end{matrix} \right]^T, \mathbf{v}_2=\left[ \begin{matrix}
	-2& 0
\end{matrix} \right]^T, \mathbf{v}_3=\left[ \begin{matrix}
	-1& 0
\end{matrix} \right]^T, \mathbf{v}_4=\left[ \begin{matrix}
	1& -4
\end{matrix} \right]^T$, we can initialize $K^{(0)}=\sum_i\mathbf{v}_i\mathbf{e}_i^T$ that makes $A-BK^{(0)}$ Schur and compartmental.
\begin{figure}[t]
    \centering
    \hfill 
    \begin{subfigure}[b]{0.49\columnwidth}
    \label{FIPM_SIPM time}
        \includegraphics[width=\linewidth]{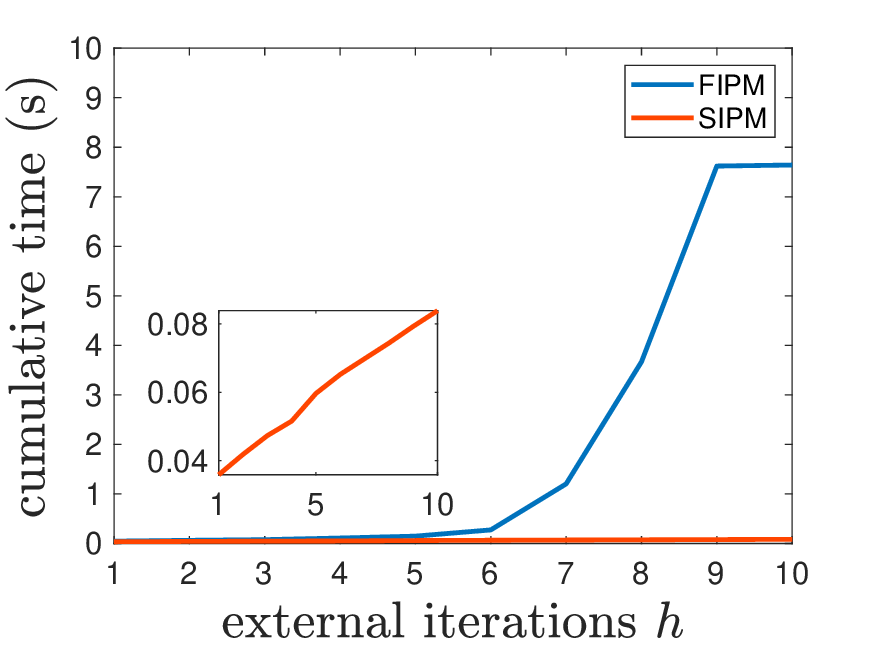}
        \caption{Cumulative running time comparison between FIPM and SIPM}
    \end{subfigure}
    \hfill
    \begin{subfigure}[b]{0.49\columnwidth}
    \label{FIPM_SIPM iter}
        \includegraphics[width=\linewidth]{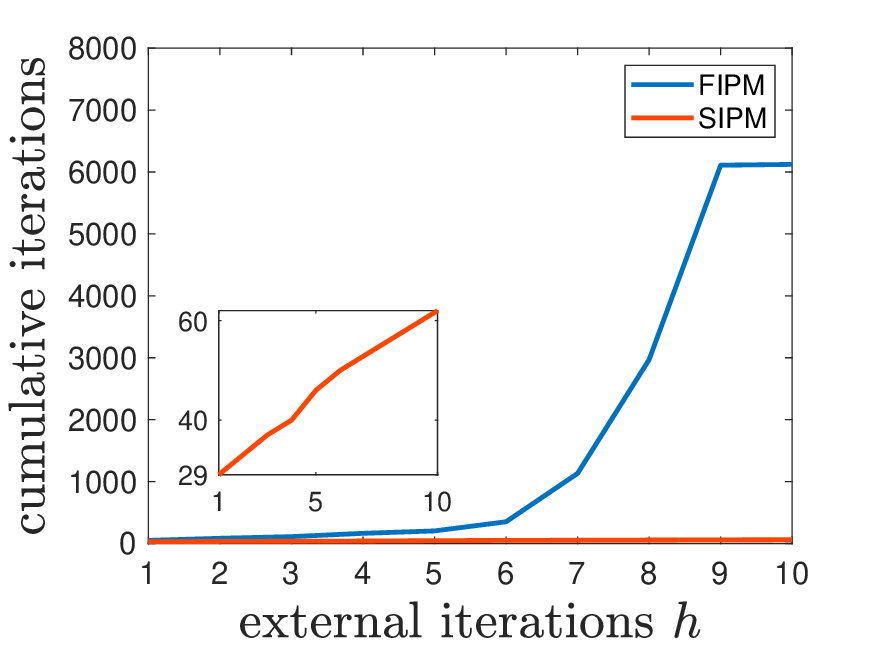}
        \caption{Cumulative iterations comparison between FIPM and SIPM}
    \end{subfigure}
    \caption{Cumulative running time and iterations comparison between FIPM and SIPM for example 1 when $N=1$.}
    \label{FIPM SIPM 1}
\end{figure}
\cref{FIPM SIPM 1}(a) shows the cumulative running time of FIPM and SIPM. After 10 external iterations, FIPM returns
\begin{equation*}
    K_F^* = \left[ \begin{matrix}
	0.6334&		0.5384&		0.6579&  0.0000\\
	0.0000&		0.5938&		0.5182&  0.5481\\
\end{matrix} \right]
\end{equation*}
with $J(K_F^*)=26.7744$, which is significantly better than $J(K^{(0)})=128.3285$. 
% Then we check
% \begin{equation*}
%     A-BK_F^*=\left[ \begin{matrix}
% 	0.4367&		0.1462&		0.0342&		0.0000\\
% 	0.1000&		0.6000&		0.0000&		0.2000\\
% 	0.4000&		0.0000&		0.8000&		0.4000\\
% 	0.0000&		0.1406&		0.0482&		0.3452\\
% \end{matrix} \right] 
% \end{equation*}
% which is Schur and compartmental.
We can easily check $A-BK_F^*$ is Schur and compartmental.
% Now we check if $K_F^*$ is a stationary point of \eqref{equivalent problem}. 
Besides, by denoting $t^{(9)}=t^*=262144$, we can check the KKT conditions hold trivially. However, as shown in \cref{FIPM SIPM 1}(a), it takes $7.2385$ seconds to converge. On the other hand, SIPM converges to 
 $K_S^*=K_F^*$ with $J(K_S^*)=26.7744$. However, it only takes $0.0862$ seconds for SIPM to converge, which is about 80 times faster than FIPM. We also compare FIPM and SIPM in terms of cumulative iterations, i.e., the cumulative number of descents. \cref{FIPM SIPM 1}(b) shows that FIPM needs $6123$ iterations while SIPM only needs $62$. 
 % Thus, in 4-dimensional systems, SIPM shows dominant advantages over FIPM. 
 % In the 
 % subsequent subsection, we will compare their performance as the system scales increase.
 
Now we compare FIPM and SIPM in various system dimensions, specifically in the context of multiple thermal systems. For $N$ thermal systems, we can write the overall system as 
\begin{equation*}
\begin{split}
    &A_N = \mathrm{blkdiag}(\overbrace{A,\ldots,A}^{N}), B_N = \mathrm{blkdiag}(\overbrace{B,\ldots,B}^{N}), \\
    &G = I_{4N\times4N},C_N = [\overbrace{C|\ldots |C}^{N}],D_N = [\overbrace{D|\ldots|D}^{N}]
\end{split}
\end{equation*}
with $K_N^{(0)}=\mathrm{blkdiag}(\overbrace{K^{(0)},\ldots,K^{(0)}}^{N})$. The simulation is shown in \cref{increasing dimension example 1}. We have also provided the simulation results of ADMM adopted in \cite{yang2024h2}. The results show that SIPM outperforms FIPM and ADMM in small-scale systems but is unsuitable for large-scale systems due to the huge computation cost of the Hessian matrix. ADMM is efficient in low-dimension systems compared to FIPM but they are comparable in high-dimension systems.
{\color{blue}
\begin{myrem}
While the convergence properties of Newton's method in general nonconvex settings remains an open problem, there are various successful applications of Newton's method in nonconvex optimization \cite{xu2020second,shi2019semismooth,ji2022globally}. Second-order methods, such as SIPM, typically achieve faster convergence by leveraging curvature information, but they require computing and inverting the Hessian matrix, which can be computationally expensive in high-dimensional systems. Therefore, when sufficient computational resources are available or the system dimension is moderate, second-order methods may be preferable for accelerating convergence. Otherwise, first-order methods, i.e., FIPM, offer a more efficient alternative.
\end{myrem}}

\begin{table*}[t]
\caption{{\color{blue}Comparison of FIPM, SIPM and ADMM with increasing system scales for example 2}}
\centering

{\color{blue}
\begin{tabular}{|l|r|r|r|r|r|r|r|}
\hline
\label{increasing dimension example 2}
 \diagbox{Method}{$N$}& 1 & 2 & 3 & 4 & 5 & 10 & 15  \\ \hline 
 FIPM& 3.8708s & 10.2208 & 17.8787s & 29.3922s & 46.9501s & 251.4954s & 573.7599s   \\ \hline
 SIPM& 0.0412s & 0.2085s & 0.6992s & 1.8664s & 3.8524s & 62.6156s & 579.0865s  \\ \hline
 ADMM \cite{yang2024h2}& 0.7187s & 6.3254s & 14.0301s & 18.5326s & 67.4557s & 727.4170s & 1842.2307s   \\ \hline
\end{tabular}}

\end{table*}
{\color{blue}
\subsection{Example 2: Leslie Matrix Model}
In this subsection, we consider the problem of pest population control, where the population dynamics are modeled by the well-known Leslie matrix model \cite{leslie1948some,liu2023proportional}. The Leslie model is widely used to describe the age-structured population growth of species, and can be expressed as \eqref{basic system} with the following system matrices
\begin{equation*}
    \begin{split}
        &A = \left[ \begin{matrix}
	b_1&		b_2&		b_3\\
	s_1&		0&		0\\
	0&		s_2&		0\\
\end{matrix} \right],B = \left[ \begin{matrix}
	0.6&		0.9\\
	0.0&		0.12\\
	0.0&		0.0\\
\end{matrix} \right], G = I_{3\times3}, \\
&C = \left[ \begin{matrix}
	I_{3\times3}\\
	0_{3\times3}\\
\end{matrix} \right], D = \left[ \begin{matrix}
	0_{4\times2}\\
	I_{2\times2}\\
\end{matrix} \right]
    \end{split}
\end{equation*}
with $C^TD=0$, $D^TD\succ0$. In this example, $x_k\in R^3$ represents the population of juvenile, immature and mature pests, respectively. $b_i$ represents the birth rate of pests in the $i$-th age group, $s_i$ represents the survival rate of pests from the $i$-th age group to the next age group. We assume that the system parameters are set as $b_1=0.25, b_2=0.6, b_3=0.56, s_1=0.35, s_2=0.25$. Our goal is to 
design a state-feedback controller to preserve the compartmental property and stability of the system, while simultaneously enhancing the system's robustness. 

\begin{figure}[t]
    \centering    \includegraphics[width=0.7\linewidth]{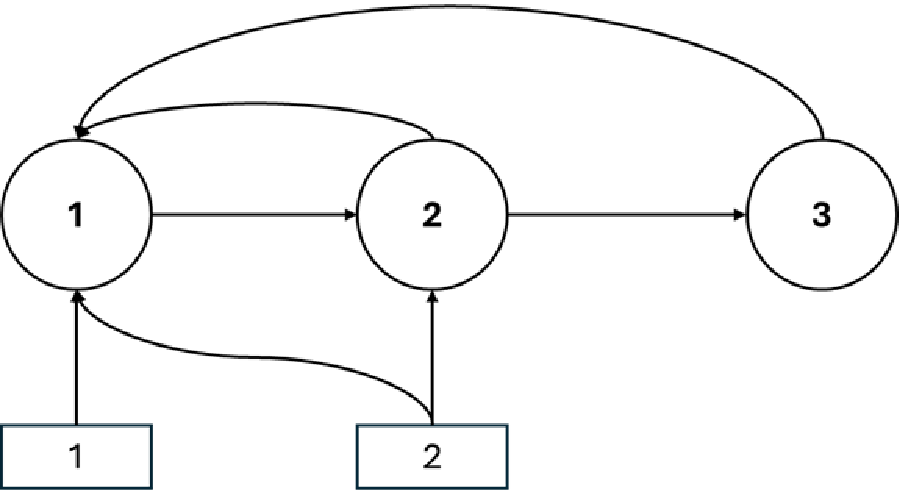}
    \caption{{\color{blue}Illustration of the Leslie model with three age groups (circle) and two inputs (rectangle)}}
    \label{Leslie model}
\end{figure}

\begin{figure}[t]
    \centering
    \hfill 
    \begin{subfigure}[b]{0.49\columnwidth}
    \label{FIPM_SIPM time}
        \includegraphics[width=\linewidth]{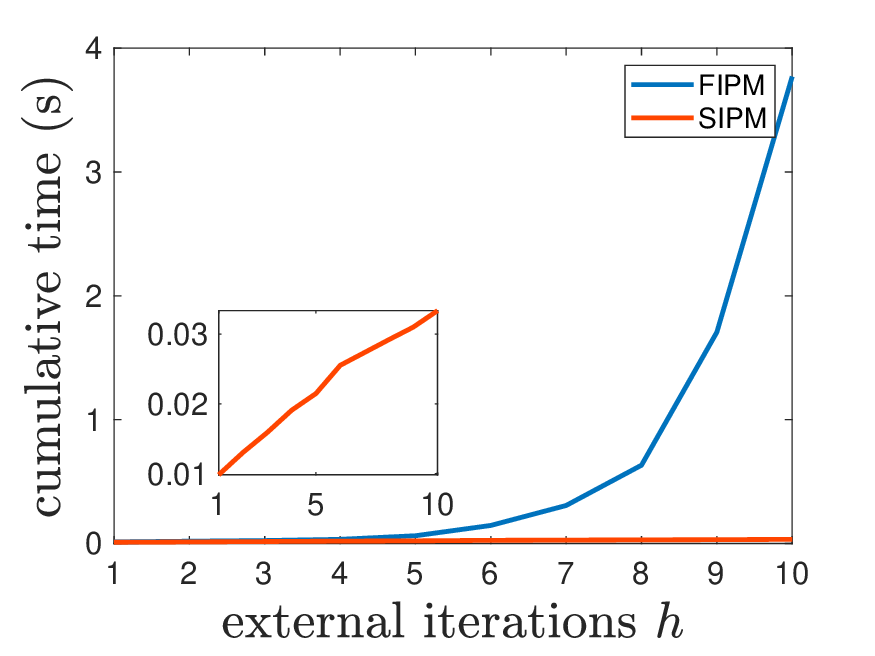}
        \caption{{\color{blue}Cumulative running time comparison between FIPM and SIPM}}
    \end{subfigure}
    \hfill
    \begin{subfigure}[b]{0.49\columnwidth}
    \label{FIPM_SIPM iter}
        \includegraphics[width=\linewidth]{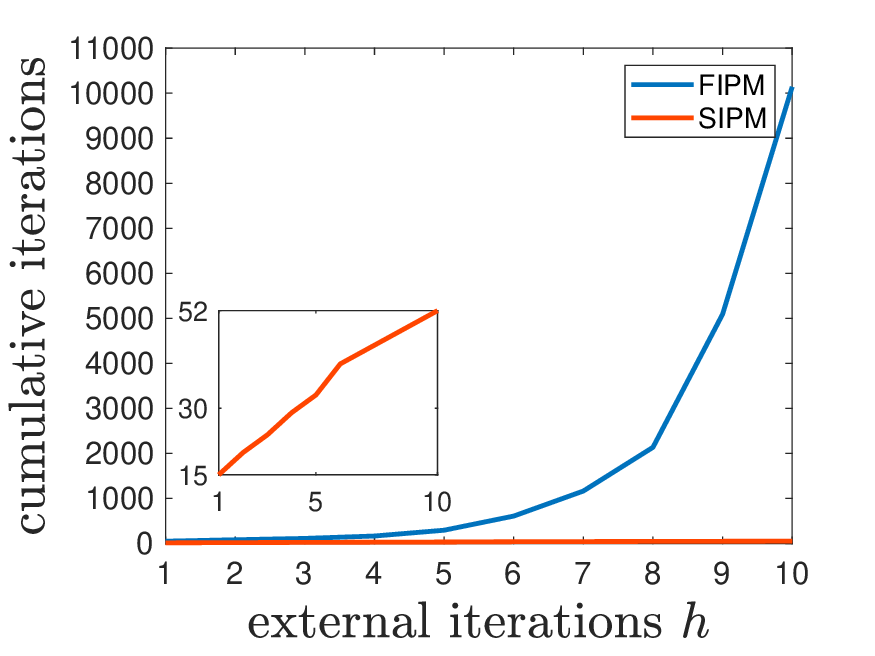}
        \caption{{\color{blue}Cumulative iterations comparison between FIPM and SIPM}}
    \end{subfigure}
    \caption{{\color{blue}Cumulative running time and iterations comparison between FIPM and SIPM for example 2 when $N=1$.}}
    \label{FIPM SIPM 2}
\end{figure}

We use the same algorithm parameters as in the previous example, i.e., $\epsilon_1=1\text{e}{-4}$, $\epsilon_2=1\text{e}{-3}$, $\epsilon_r=1\text{e}{-9}$, $\delta=1$, and $\mu = 4$ for both FIPM and SIPM. It is straightforward to verify that $A-I$ is irreducible, which allows us to apply Lemma 7 from \cite{valcher2018state} to construct an initial controller $K^{(0)}$ such that $A-BK^{(0)}$ is both Schur and compartmental. Specifically, by selecting $\mathbf{v}_1=\left[ \begin{matrix} 1 & -1 \end{matrix} \right]^T$ and $\mathbf{v}_3=\left[ \begin{matrix} -2 & 1 \end{matrix} \right]^T$, we set $K^{(0)}=\frac{1}{2}\mathbf{v}_1\mathbf{e}_1^T+\frac{1}{2}\mathbf{v}_3\mathbf{e}_3^T$. 
After 10 external iterations, FIPM yields the following controller:
\begin{equation*}
    K_F^* = \left[ \begin{matrix}
    0.0518 & 0.3055 & 0.2804 \\
    0.1856 & 0      & 0      \\
\end{matrix} \right]
\end{equation*}
with a corresponding performance value $J(K_F^*)=3.8429$, which is a significant improvement over the initial value $J(K^{(0)})=5.7594$. It can be readily verified that $A-BK_F^*$ is both Schur and compartmental, and the KKT conditions are satisfied. The total computation time for FIPM in this case is $3.8708$ seconds, as shown in \cref{FIPM SIPM 2}(a). 
In contrast, SIPM converges to the same controller $K_S^*=K_F^*$ with only $0.0412$ seconds, demonstrating a remarkable speedup of approximately 94 times compared to FIPM. The performance value for SIPM is also $J(K_S^*)=3.8429$. Moreover, we oberve that SIPM requires only $52$ iterations to converge, while FIPM requires $10141$ iterations, as shown in \cref{FIPM SIPM 2}(b).

We further extend the simulation to larger system scales by considering multiple coupled Leslie models. For $N$ compartments, the overall system matrices are constructed as follows:
\begin{equation*}
\begin{split}
    &A_N = \mathrm{blkdiag}(\overbrace{A,\ldots,A}^{N}), \quad B_N = \mathrm{blkdiag}(\overbrace{B,\ldots,B}^{N}), \\
    &G_N = I_{3N\times3N}, \quad C_N = [\overbrace{C|\ldots|C}^{N}], \quad D_N = [\overbrace{D|\ldots|D}^{N}]
\end{split}
\end{equation*}
with the initial controller $K_N^{(0)} = \mathrm{blkdiag}(\overbrace{K^{(0)},\ldots,K^{(0)}}^{N})$.

We assess the computational performance of FIPM, SIPM, and ADMM for $N = 1, 2, 3, 4, 5, 10, 15$ compartments. The results are summarized in \cref{increasing dimension example 2}. SIPM consistently achieves faster convergence than both FIPM and ADMM for small and moderate system sizes. However, as the system scale increases, the computational cost of SIPM grows rapidly due to the complexity of Hessian matrix calculations, making it less suitable for very large-scale problems. In contrast, FIPM exhibits more moderate growth in computation time and remain practical for larger systems. In summary, SIPM is highly efficient for small to medium-scale compartmental systems, while FIPM is preferable for large-scale applications.
}

\section{Conclusion}\label{Conclusion}
In this paper, we studied the $H_2$ optimal control for compartmental systems. We proposed a novel problem transformation and established an equivalent new optimization problem with closed and polyhedron constraints. We provided FIPM and a novel SIPM to solve the problem, both are guaranteed to converge to a stationary point of the new problem. Meanwhile, we propose an initialization method to guarantee the interior property of initial values. Finally, we conducted thorough simulations to compare FIPM, SIPM, and ADMM.

% {\color{blue}and compare them with existing ADMM}. 
% The technique and analysis c
% an be easily extended to 
% static output-feedback 
% control, which is left for future research.
\footnotesize{
\bibliographystyle{IEEEtran}
\bibliography{reference}
}
\end{document}